\documentclass[11pt]{article}
\usepackage{amsmath, amsfonts, amsthm, amssymb}

\newcommand{\ignore}[1]{}

\def\@begintheorem#1#2{\par\bgroup{\sc #1\ #2. }\it\ignorespaces}
\def\@opargbegintheorem#1#2#3{\par\bgroup{\sc #1\ #2\ (#3). }
\it\ignorespaces}
\def\@endtheorem{\egroup}
\newtheorem{theorem}{Theorem}[section]
\newtheorem{corollary}[theorem]{Corollary}
\newtheorem{lemma}[theorem]{Lemma}
\newtheorem{proposition}[theorem]{Proposition}
\newtheorem{problem}[theorem]{Problem}
\newtheorem{example}[theorem]{Example}
\newtheorem{algorithm}[theorem]{Algorithm}
\newtheorem{definition}[theorem]{Definition}
\newcommand{\bt}[1]{\begin{theorem}\label{#1}}
\newcommand{\bc}[1]{\begin{corollary}\label{#1}}
\newcommand{\bl}[1]{\begin{lemma}\label{#1}}
\newcommand{\bp}[1]{\begin{proposition}\label{#1}}
\newcommand{\bpro}[1]{\begin{problem}\label{#1}}
\newcommand{\be}[1]{\begin{example}\rm\label{#1}}
\newcommand{\ba}[1]{\begin{algorithm}\rm\label{#1}}
\newcommand{\bd}[1]{\begin{definition}\rm\label{#1}}
\newcommand{\bpr}{\begin{proof}}
\newcommand{\et}{\end{theorem}}
\newcommand{\ec}{\end{corollary}}
\newcommand{\el}{\end{lemma}}
\newcommand{\ep}{\end{proposition}}
\newcommand{\epro}{\end{problem}}
\newcommand{\ee}{\end{example}}
\newcommand{\ea}{\end{algorithm}}
\newcommand{\ed}{\end{definition}}
\newcommand{\epr}{\end{proof}}

\def\Z{\mathbb{Z}}
\def \supp {{\rm supp}}
\def \fr {{\rm F}}

\begin{document}

\title{\bf Intractability of approximate multi-dimensional
nonlinear optimization on independence systems}
\author{Jon Lee \and Shmuel Onn \and Robert Weismantel}
\date{}
\maketitle

\begin{abstract}
We consider optimization of nonlinear objective functions
that balance $d$ linear criteria over $n$-element independence
systems presented by linear-optimization oracles.
For $d=1$, we have previously shown that an $r$-best approximate
solution can be found in polynomial time. Here, using an extended
Erd\H{o}s-Ko-Rado theorem of Frankl, we show that for $d=2$,
finding a $\rho n$-best solution requires exponential time.
\end{abstract}

\section{Introduction}

Given system $S\subseteq\{0,1\}^n$, integer $d\times n$
matrix $W$, and function $f:\Z^d\rightarrow\Z$,
consider the problem of minimizing the nonlinear
composite function $f(Wx)$ over $S$, that is,
\begin{equation}\label{NCO}
\min\{f(Wx)\ :\ x\in S\}\ .
\end{equation}
This problem can be interpreted as multi-criteria optimization,
where row $W_i$ of $W$ gives a linear function $W_ix$ representing
the value of feasible point $x\in S$ under criterion $i$,
and the objective value $f(Wx)=f(W_1x,\dots,W_dx)$
is the balancing of these $d$ criteria.

Assume we can do linear optimization over
$S$ to begin with, namely $S$ is presented by a
{\em linear-optimization oracle}, which queried on
$w\in \Z^n$, solves $\max\{wx\,:\,x\in S\}$.
For restricted systems $S$, such as matroids and
matroid intersections, or restricted functions $f$,
such as concave functions, problem (\ref{NCO})
can be solved in polynomial time \cite{BLMORWW,BLOW}.
A comprehensive description of the state of the
art on this area can be found in \cite{Onn}.

Here we continue our investigation from \cite{LOW}
of problem (\ref{NCO}) where $S$ is an arbitrary
{\em independence system}, that is, $S$ nonempty,
and $x\leq y\in S$ with $x\in\{0,1\}^n$ imply $x\in S$.

A feasible point $x^*\in S$ is called
an {\em $r$-best solution} of problem (\ref{NCO})
provided there are at most $r$ better objective function
values attainable by other feasible points, that is,
$$\left|\left\{f(Wx)\ :\ f(Wx)<f(Wx^*)
\,,\ x\in S\right\}\right|\ \leq\ r\ .$$
So it provides a suitable approximation to (\ref{NCO}).
In particular, a $0$-best solution is optimal.

In \cite{LOW}, the case of $d=1$ was considered,
that is, the problem $\min\{f(wx):x\in S\}$ with $w\in\Z^n$.
It was shown that for any fixed positive integers $a_1,\dots,a_p$
there is a polynomial time algorithm that, given any
$w\in\{a_1,\dots,a_p\}^n$, provides an $r(a_1,\dots,a_p)$-best
solution to the problem, where $r(a_1,\dots,a_p)$ is a constant
related to Frobenius numbers of some of the $a_i$. In particular,
for any $p=2$ integers, $r(a_1,a_2)=\fr(a)$ is the Frobenius number.

In this note we consider the problem in dimension $d=2$.
We restrict attention to $2\times n$ matrices $W$ which are
$\{0,1\}$-valued. Then the {\em image} of $S$ under $W$ satisfies
\begin{equation}\label{set}
WS\ :=\ \{Wx\, :\, x\in S\}\ \subseteq\ \{0,1,\dots,n\}^2\ .
\end{equation}
Therefore, the problem of computing the optimal
objective function {\em value} of (\ref{NCO}) is seemingly reducible
to computing the image $WS$ by checking if $y\in WS$ for each
of the $(n+1)^2$ points $y$ in the set on right-hand side
of (\ref{set}) and determining the minimum value of $f$ over $WS$.
Unfortunately, this so called {\em fiber problem}, of checking
if $y\in WS$, is computationally hard. In particular, already
for $S$ the set of (indicators of) matchings in a bipartite graph,
over which linear optimization is easy, this problem includes
as a special case the notorious {\em exact matching problem}
whose complexity is long open \cite{PY}.

Here we show that there is a universal positive constant $\rho$
such that, already for $d=2$, matrix $W$ each column of which is
one of the two unit vectors in $\Z^2$, and very simple explicit
function $f$ supported on $\{0,1,\dots,n\}^2$, there is no
polynomial time algorithm that can produce even a $\rho n$-best
solution of problem (\ref{NCO}) for every independence
system $S\subseteq\{0,1\}^n$, let alone find a constant
$r$-best or optimal solution. Our construction
makes use of a beautiful extension of the classical
Erd\H{o}s-Ko-Rado theorem due to Frankl \cite{Fra}.

It is interesting whether our construction could be refined
to shed some light on the exact matching and related open problems
of \cite{PY}, and whether other natural oracles for $S$ could
lead to polynomial time solution of problem (\ref{NCO})
in dimensions $d=2$ and higher.

\section{A $\rho n$-best solution cannot be found in polynomial time}
\label{ApproximationLowerBound}

\bt{exponential_lower_bound}
There exists a universal positive constant $\rho$ such that no
polynomial time algorithm can compute a $\rho n$-best solution
of the $2$-dimensional nonlinear optimization problem
$\min\{f(Wx)\,:\,x\in S\}$ over every independence system
$S\subseteq\{0,1\}^n$ presented by a linear-optimization
oracle, with $W$ an integer $2\times n$ weight matrix each
column of which is one of the unit vectors in $\Z^2$, and
$f$ an explicit function supported on $\{0,1,\dots,n\}^2$.

\vskip.2cm\noindent
In fact, the following explicit statement holds.
Let $l$ be any positive integer with $l\geq 2^{10}$,
$k:=7l$, $m:=8l^2$, $n:=2m$, and $\rho:={1\over 17}$.
Let $W$ be the $2\times n$ matrix with first $m$ columns
the unit vector ${\bf 1}_1$ and last $m$ columns
the unit vector ${\bf 1}_2$.
Define $f$ on $\Z^2$ explicitly by
\begin{equation}\label{function}
f(y)\ =\ f(y_1,y_2)\ :=\
\left\{
  \begin{array}{ll}
    (y_1-k)-l(y_2-k)-1 & \mbox{if}\ k+1\leq y_1,y_2\leq k+l\ ,\\
     0 & \mbox{otherwise}\ .
  \end{array}
\right.
\end{equation}
Then at least $2^{{1\over4}\sqrt n}$ queries to the oracle
of $S$ are needed to compute a ${1\over 17}n$-best solution.
\et
\bpr
Let $l\geq 2^{10}$ be a positive integer, $k,m,n,\rho$ and
$W$ as above, and $f$ as in (\ref{function}) above.
It is more convenient here to work with set systems over ground set
$N:=\{1,\dots,n\}$ rather than sets of vectors in $\{0,1\}^n$.
As usual, vectors $x\in\{0,1\}^n$ are in bijection with subsets
$X\subseteq N$ with corresponding elements satisfying $X=\supp(x)$
the support of $x$ and $x={\bf 1}_X$ the indicator of $X$.
So we replace each $S\subseteq \{0,1\}^n$
by the set system ${\cal S}:=\{X=\supp(x)\,:\,x\in S\}$.
Also, for $c\in\Z^n$ and $X\subseteq N$ we write $cX:=c{\bf 1}_X$.
Let $N_1\uplus N_2=N$ be the natural equipartition of the
ground set defined by
$N_1:=\{1,\dots,m\}$ and $N_2:=\{m+1,\dots,2m\}$.
For each subset $X\subseteq N$ of the ground set we write
$X_1:=X\cap N_1$, $X_2:=X\cap N_2$, with $X=X_1\uplus X_2$ the
naturally induced partition of $X$.

The image of $X=X_1\uplus X_2$ is denoted by $WX:=W{\bf 1}_X$ and is
equal to $(|X_1|,|X_2|)$. The image of a set system $\cal S$ over
$N$ is $W{\cal S}:=\{WX:X\in{\cal S}\}$.
We use several set systems over $N$, defined as follows.
First, for each pair of integers $0\leq y_1,y_2\leq m$, let
$${\cal S}_{y_1,y_2}\ :=\
\left\{X=X_1\uplus X_2\ :\ |X_1|=y_1\,,\ |X_2|=y_2\right\}\ .$$
Next, let
$${\cal S}^*\ :=\
\left\{X\ : \ (|X_1|,|X_2|)\leq (m,k)
\ \ \mbox{or}\ \ (|X_1|,|X_2|)\leq (k,m)\right\}\ .$$
Then $\cal S^*$ is an independence system whose image is given by
\begin{eqnarray*}
W{\cal S}^* & = &
\left\{(y_1,y_2)\in\Z^2_+\ :\ (y_1,y_2)\leq (m,k)\ \
\mbox{or}\ \ (y_1,y_2)\leq (k,m)\right\}\ .
\end{eqnarray*}
Moreover, the objective function value of every $X\in\cal S^*$,
and hence in particular of every $\rho n$-best solution of the
minimization problem over $\cal S^*$, satisfies $f(WX)=0$.

Next, for each $Y\in{\cal S}_{k+l,k+l}$, let
$${\cal S}_Y\ :=\
{\cal S}^*\ \cup\ \left\{X\ : \ X\subseteq Y\right\}\ .$$
Then ${\cal S}_Y$ is also an independence system, with image
\begin{eqnarray*}
W{\cal S}_Y & = & W{\cal S}^*\
\uplus\ \left\{(y_1,y_2)\ :\
(k+1,k+1)\leq(y_1,y_2)\leq(k+l,k+l)\right\}\ .
\end{eqnarray*}
Moreover, the objective function values of the points in
${\cal S}_Y\setminus{\cal S}^*$, whose images lie in
$W{\cal S}_Y\setminus W{\cal S}^*$, attain exactly all
$l^2={1\over 16}n>\rho n$ values $-1,-2,\dots,-l^2$,
and so the value of every
$\rho n$-best solution of the minimization problem
over ${\cal S}_Y$ satisfies $f(WX)\leq -1$.

\vskip.3cm
For each vector $c\in\Z^n$ and each pair $1\leq i_1,i_2\leq l$, let
$${\cal T}_{i_1,i_2}(c)\ :=\ \left\{Z\in{\cal S}_{k+i_1,k+i_2}
\ :\ cZ>\max\{cX\ :\ X\in{\cal S}^*\}\right\}\ .$$

\noindent \emph{Claim:} For every $c\in\Z^n$ and every
pair $1\leq i_1,i_2\leq l$, we have
$${{|{\cal T}_{i_1,i_2}(c)|}}\ \leq \ {m\choose l}{m\choose k+l}\ .$$

\noindent \emph{Proof of Claim:}
Consider any pair
$U=U_1\uplus U_2\,,V=V_1\uplus V_2\in {\cal T}_{i_1,i_2}(c)$.
We now show that either
$|U_1\cap V_1|\geq k+1$ or $|U_2\cap V_2|\geq k+1$.
Suppose, indirectly, this is not so. Put
\begin{eqnarray*}
X &:=& (U_1\cap V_1)\uplus(U_2\cup V_2), \\
Y &:=& (U_1\cup V_1)\uplus(U_2\cap V_2).
\end{eqnarray*}
Then $|U_1\cap V_1|\leq k$ and $|U_2\cup V_2|\leq m$
imply $X\in{\cal S}^*$, and $|U_1\cup V_1|\leq m$ and
$|U_2\cap V_2|\leq k$ imply $Y\in{\cal S}^*$.
We then obtain the following contradiction,
$$0\ <\ cU-cX\ =\
c(U_1\setminus V_1)-c(V_2\setminus U_2)\ =\ cY-cV\ <\ 0 .$$
So indeed, for every pair
$U=U_1\uplus U_2\,,V=V_1\uplus V_2\in
{\cal T}_{i_1,i_2}(c)\subseteq{\cal S}_{k+i_1,k+i_2}$,
either $|U_1\cap V_1|\geq k+1$ or $|U_2\cap V_2|\geq k+1$.
Therefore, we can now apply the extended Erd\H{o}s-Ko-Rado theorem
for direct products of Frankl \cite[Theorem 2]{Fra}, which implies
$$\hskip-0.25cm
{{|{\cal T}_{i_1,i_2}(c)|}\over{|{\cal S}_{k+i_1,k+i_2}|}}\ \leq \
\max\left\{{m-(k+1)\choose (k+i_1)-(k+1)}\left/{m\choose k+i_1}\right.,
{m-(k+1)\choose (k+i_2)-(k+1)}\left/{m\choose k+i_2}\right.\right\}
$$
from which it is easy to conclude that, as claimed,
$${{|{\cal T}_{i_1,i_2}(c)|}}\ \leq\ {m\choose l}{m\choose k+l}\ .$$

\vskip.4cm
We continue with the proof of our theorem.
Since $k=7l$, $m=8l^2$ and $l\geq 2$ we get
\begin{eqnarray*}
{m\choose k+l}\left/{m\choose l}^3\right.\ \ =\ \
{8l^2\choose 8l}\left/{8l^2\choose l}^3\right. \ \ \geq\ \
\left({{4l^2}\over 8l}\right)^{8l}\left/(8l^2)^{3l}\right.
\ \ \geq\ \ (2^{-9}l)^{2l}\ \ .
\end{eqnarray*}
Therefore
$$|{\cal S}_{k+l,k+l}|\ \ =\ \ {m\choose k+l}{m\choose k+l}
\ \ \geq\ \ (2^{-9}l)^{2l} {m\choose l}^3{m\choose k+l}\ \ .$$

Consider any algorithm attempting to obtain a $\rho n$-best solution
to the nonlinear optimization problem over any system $\cal S$,
and let $c^1,\dots,c^q\in\Z^n$ be the sequence of queries
to the oracle of $\cal S$ made by the algorithm.
For each pair $1\leq i_1,i_2\leq l$ and each $Z\in {\cal T}_{i_1,i_2}(c^p)$, the number
of $Y\in {\cal S}_{k+l,k+l}$ containing $Z$, and hence satisfying
$Z\in{\cal S}_Y$, is
$${m-(k+i_1)\choose l-{i_1}}{m-(k+i_2)\choose l-{i_2}}
\ \ \leq \ \ {m\choose l}^2\ \ .$$
So the number of $Y\in {\cal S}_{k+l,k+l}$ containing some $Z$
which lies in some ${\cal T}_{i_1,i_2}(c^p)$ is at most
$$\sum_{p=1}^q\sum_{i_1=1}^l\sum_{i_2=1}^l
{m\choose l}^2|{\cal T}_{i_1,i_2}(c^p)|
\ \ \leq\ \ q l^2 {m\choose l}^3{m\choose k+l}\ .$$
Therefore, if the number of oracle queries satisfies
$q<l^{-2}(2^{-9}l)^{2l}$, then there exists some
$Y\in{\cal S}_{k+l,k+l}$ which does not contain any $Z$ in any
${\cal T}_{i_1,i_2}(c^p)$. This means that any $Z\in{\cal S}_Y$
satisfies $c^p Z\leq\max\{c^p X:X\in{\cal S}^*\}$. Hence, whether
the linear-optimization oracle presents $\cal S^*$ or ${\cal S}_Y$,
on each query $c^p$ it can reply with some $X^p\in{\cal S}^*$
attaining
$$c^p X^p\ =\ \max\{c^p X\ :\ X\in {\cal S}^*\}
\ =\ \max\{c^p X\ :\ X\in{\cal S}_Y\}\ .$$
So the algorithm cannot tell whether the oracle presents
$\cal S^*$ or ${\cal S}_Y$, whether the image is $W{\cal S}^*$
or $W{\cal S}_Y$, and whether the objective function value of
every $\rho n$-best solution is zero or negative, let alone
compute any $\rho n$-best solution. Therefore, with $l\geq 2^{10}$,
every algorithm which can produce a $\rho n$-best solution for
the $2$-dimensional nonlinear optimization problem (\ref{NCO})
over every system $\cal S$ must make at least an exponential number
$$q\ \ \geq\ \ l^{-2}(2^{-9}l)^{2l}\ \ \geq\ \ l^{-2} 2^{2l}
\ \ >\ \ 2^l\ \ =\ \  2^{{1\over4}\sqrt n}$$ of queries to
the oracle presenting $\cal S$ and
therefore cannot run in polynomial time.
\epr

\vskip.6cm\noindent {\small Jon Lee}\newline
\emph{IBM T.J. Watson Research Center, Yorktown Heights, USA}

\vskip.3cm\noindent {\small Shmuel Onn}\newline
\emph{Technion - Israel Institute of Technology, Haifa, Israel}

\vskip.3cm\noindent {\small Robert Weismantel} \newline
\emph{ETH, Z\"urich, Switzerland}

\end{document}